\documentclass{amsart}
\usepackage{graphics} 
\usepackage{psfrag} 
\usepackage{graphicx}
\newtheorem{lemma}{Lemma}

\def\R{\mathbb{R}}
\def\N{\mathbb{N}}

\newcommand{\h}{\mathcal{H}}
\def\guilo{«~}
   \def\guilf{~»}
    \def\andname{, }
    
\def\volumename{vol.~}
 \def\numbername{\numero}
    \def\numbername{\numero}
  \def\numero{n\@umer{o}}%
\def\Numero{N\@umer{o}}%
\title{ An intermediate targets method for time parallelization in optimal control}

\author[Y. Maday]{Yvon Maday}
\address{ Universit\'e Pierre et Marie Curie-Paris6, UMR
7598, Laboratoire J.-L. Lions, Paris, F-75005 France, and  Division of Applied Mathematics, Brown University
182 George Street, Providence, RI 02912, USA,}
\email{maday@ann.jussieu.fr}
\author[J. Salomon]{Julien Salomon}
\address{Universit\'e Paris-Dauphine, UMR 7534, CEREMADE, Paris, F-75016 France,}
\email{salomon@ceremade.dauphine.fr}
\author[K. Riahi]{Kamel Riahi}
\address{Universit\'e Pierre et Marie Curie-Paris6, UMR
7598, Laboratoire J.-L. Lions, Paris, F-75005 France.}
\email{riahi@ann.jussieu.fr}

\begin{document}
\centerline{}
\maketitle
\begin{abstract}
In this paper, we present a method that enables to solve in parallel the
  Euler-Lagrange system associated with the optimal control of a parabolic equation. 
Our approach is based on an
iterative update of a sequence of intermediate targets and
gives rise independent sub-problems that can be solved in parallel. 
Numerical experiments show the efficiency of our method.


Dans cet article, on pr\'{e}sente une m\'ethode permettant une
  parall\'elisation en temps de la r\'{e}solution des \'equations
  d'Euler-Lagrange associ\'{e}es \`{a} un probl\`{e}me de
  contr\^{o}le optimal dans le cas parabolique. 
Notre approche est  bas\'{e}e sur une mise \`{a} jour it\'{e}rative de cibles
  interm\'{e}diaires et donne lieu \`{a} des sous-probl\`{e}mes de
  contr\^{o}le ind\'{e}pendants. Les r\'{e}sultats num\'{e}riques prouvent l'efficacit\'{e} de la m\'{e}thode.

\end{abstract}



\section{Introduction:}
In the last decade, time domain decomposition has been exploited to
accelerate the simulation of systems ruled by time dependent
partial differential equations~\cite{L}.  Among others,
the parareal algorithm~\cite{LMT} or multi-shooting
schemes~\cite{BZ} have shown excellent results. In the framework of
optimal control, this approach has been used to control parabolic
 systems~\cite{MT}, ~\cite{MSS}.  In this paper, we introduce a new approach
 to tackle such problems. The strategy we follow is based on the concept of
 target trajectory that has been introduced in the case of hyperbolic
 systems in~\cite{MST}. Because of the irreversibility of parabolic
 equations, a new definition of this trajectory is considered.
 It enables us to define at each bound of the
 time sub-domains relevant initial conditions and intermediate
 targets, so that the initial problem is split into independent
 optimization problems.\\

We now introduce some notations. Given $d\in\N$, we consider the optimal control problem associated
with a heat equation defined on a compact set $\Omega\subset\R^d$ and
a time interval interval $I=[0,T]$. The control applies on a subset
$\Omega_c\subset\Omega$ and belongs to the Hilbert space
$\h:=L^{2}(I;L^{2}(\Omega_c))$ whose scalar product and norm are
denoted  by $\langle.,.\rangle_{\h}$ and $\|.\|_{\h}$ respectively. We
also denote by $\|.\|$ and $\|.\|_c$ the 
norms associated with the spaces $L^{2}(\Omega)$ and  $L^{2}(\Omega_c)$
respectively. Given a function $\varphi$
defined on $I\times\Omega$, we denote its restriction to a sub-domain $I'\times\Omega
\subset I\times\Omega$ by $\varphi_{|_{I'}}$.\\
The paper is organized as follows.
The optimal control problem is presented in Sec.~\ref{sec:part2} ; Sec.~\ref{sec:part3}
 is devoted to the description of the time
parallelization setting. The algorithm is given in
Sec.~\ref{sec:part4}. Finally, numerical results are presented in Sec.~\ref{sec:num}.
\section{Optimal control problem}\label{sec:part2}
Given $\alpha>0$, $T>0$, a target state $y_{target}\in L^{2}(\Omega)$, we consider the optimal control
problem : Find   $v^\star\in\h$
\begin{equation*}
v^\star:=\hbox{argmin}_{v\in\h} J(v),
\end{equation*}
where $J$ is the quadratic cost functional defined by
\begin{equation}\label{reference-functional}
 J(v):=\frac{1}{2}\|y(T)-y_{target}\|^2 + \frac{\alpha}{2}\| v \|_{\h}^2.
\end{equation}
The state variable $y$ depends linearly on the control $v$ through the evolution equation
\begin{equation}\label{heat}
\partial_t {y} - \nu \Delta y =Bv,
\end{equation} 
with an initial data $y_{0}\in L^2(\Omega)$. In this equation
 $\Delta$ denotes the Laplace operator, $\nu$ is the diffusion coefficient and $B$ is the
injection from $L^{2}(\Omega_c)$ to $L^{2}(\Omega)$, so that
$Bv\in L^{2}(I;L^{2}(\Omega))$.\\ 
The corresponding optimality system reads
\begin{align}
&\left\{\begin{array}{ccl}\label{os1}
\partial_t{y} - \nu \Delta y&=& Bv  \qquad \hbox{on} \ I\times\Omega\\
y(0)&=&y_0,
\end{array}
\right.\\
&\left\{\begin{array}{ccl}\label{os2}
\partial_t{p} + \nu \Delta p &=& 0 \qquad \hbox{on} \ I\times\Omega\\
p(T) &=&  y(T)-y_{target},
\end{array}
\right.\\
&\label{os3}\alpha v + B^*p =0,
\end{align}
where $B^*$ is the adjoint operator of $B$.\\
Note that for any $\alpha>0$, the functional $J$ is strictly convex,
so that the system (\ref{os1}--\ref{os3}) has a unique
solution $v^\star$. We denote by $y^\star$, $p^\star$ the
associated state and adjoint state.  
\section{Time parallelization setting}\label{sec:part3}
We now aim at solving in parallel the coupled system corresponding to
Equations (\ref{os1}--\ref{os3}). To do this, we decompose the
interval $I$ into subintervals and introduce a set of intermediate
target states so that the initial problem is replaced by a set of
smaller independent optimal control problems. The resolution of these problems is  
achieved using an inner loop, whereas the intermediate states are
updated by an outer loop.\\
We start with the definition of the target states. Given a control $v\in\h$ and its corresponding
state $y$ and adjoint state $p$, we define the {\it target trajectory} by:
\begin{equation}\label{target_trajectory}
\chi(v) = y(v) - p(v) \qquad \hbox{on}\quad I\times\Omega
\end{equation}
In what follows and for the sake of simplicity, we omit the dependence of $y$, $p$ and $\chi$ on the control $v$ in the notations. The introduction of this trajectory is motivated by the following result.
\begin{lemma}\label{lem1}We keep the previous notations. Let $\tau\in]0,T[$, and
      the optimal control problem: Find $w_{\tau}^{\star}\in\h$ such
      that 
\begin{equation*}
w^{\star}_\tau:=argmin_{w\in\h}{J}_{\tau}(w),
\end{equation*}
where $${J}_{\tau}(w):=\dfrac{1}{2}\|y(\tau)-\chi^{\star}(\tau)\|^{2} +\dfrac{\alpha}{2}\|w\|^{2}_{L^{2}([0,\tau];L^{2}(\Omega_c))}$$
with $y(\tau)$ the solution of the Equation \eqref{heat}. We have
 $$w_{\tau}^{\star}=v^{\star}_{|_{[0,\tau]}}$$
\end{lemma}

{\bf Proof: } Thanks to the uniqueness of the solution of the optimization problem associated to ${J}_{\tau}$ , it is
sufficient to show that $v^{\star}_{|_{[0,\tau]}}$ is a solution of its optimality
system. This one is given by Equation~\eqref{os1} (with $v=w_\tau$) restricted to $[0,\tau]\times\Omega$ and 

\begin{equation}\label{pos2bis}
\left\{\begin{array}{ccl}
\partial_t{\tilde p} + \nu \Delta \tilde p &=& 0 \qquad \hbox{on} \ [0,\tau]\times\Omega\\
\tilde p(\tau) &=&  y(\tau)-\chi(\tau),
\end{array}
\right.
\end{equation}

\begin{equation}\label{pos3bis}
\alpha w_\tau + B^*\tilde p =0.
\end{equation}

First, note that $y_{|_{[0,\tau]}}^\star$ obviously satisfies Equation~\eqref{os1} restricted to $[0,\tau]\times\Omega$
with $v=v^\star_{|_{[0,\tau]}}$. It directly follows from the definition of $\chi^\star$
(see \eqref{target_trajectory}), that:
\begin{equation*}
p^\star(\tau) =  y^{\star}(\tau)-\chi^\star(\tau),
\end{equation*}
so that $ p^\star_{|_{[0,\tau]}}$ satisfies~\eqref{pos2bis}.
Finally,~Equation \eqref{pos3bis} is a consequence of~\eqref{os3}. The result follows.
$\hfill \square$

Given $N\geq 1$, we decompose the interval $I=[0,T]$ into subintervals
$I=\cup_{n=0}^{N-1}I_n$, where $I_n=[t_n,t_{n+1}]$,
$t_0=0<t_1<...<t_{N-1}<t_N=T$. We also introduce the spaces
$\h_n:=L^{2}(I_n;L^{2}(\Omega_c))$ and the corresponding scalar 
product $\langle.,.\rangle_{\h_n}$ and norm $\|.\|_{\h_n}$. In this
framework, given $v_n\in\h_n$ we define $ v^\star_{n}$ as follows
\begin{equation}\label{min-loc}
v_n^\star:=\hbox{argmin}_{v_n\in \h_n}J^v_n(v_n),
\end{equation}
with
\begin{equation}\label{subpb}
 J^v_n(v_n):=\frac12 \|y_n(t_{n+1})-\chi(t_{n+1})\|^2 +
\frac{\alpha}{2}\|v_n\|_{\h_n}^2,
\end{equation}
where $\chi$ is associated to $v$ through the definition~\eqref{target_trajectory}.
In this functional, the state $y_n$ is defined by
\begin{equation}\label{pos1}
\left\{\begin{array}{ccl}
\partial_t{y_n} - \nu \Delta y_n&=& Bv_n  \qquad \hbox{on} \quad I_n\times\Omega\\y_n(t_n)&=&y(t_n).
\end{array}
\right.
\end{equation}
These subproblems have the same structure as the original one and are
also strictly convex. Note also that their definitions depend on the
control $v$ through the target trajectory, hence the notation $J_n^v$.\\
The optimality system associated with these
minimization problems are given by Equation~\eqref{pos1} and 

\begin{equation}\label{pos2}
\left\{\begin{array}{ccl}
\partial_t{p_n} + \nu \Delta p_n &=& 0 \qquad \hbox{on} \ I_n\times\Omega\\
p_n(t_{n+1}) &=&  y(t_{n+1})-\chi(t_{n+1}),
\end{array}
\right.
\end{equation}

\begin{equation}\label{pos3}
\alpha v_n + B^*p_n =0.
\end{equation}
\begin{lemma}\label{lem2}
We keep the previous notations. Denote by $\chi^\star$  the target trajectory defined by
Equation~\eqref{target_trajectory} with $y=y^\star$ and $p=p^\star$
and  by $y_n^\star,p_n^\star,v_n^\star$ the solutions of Equations
(\ref{pos1}--\ref{pos3}) associated with $v^\star$. One has:
$$ v_n^\star = v^\star_{|_{I_n}}.$$
\end{lemma}
The proof of this result follows the lines of Lemma~\ref{lem1} and is left as an
exercise to reader.\\

\section{Algorithm}\label{sec:part4}
We are now in the position to propose a time parallelized procedure
to solve Equations \eqref{os1}--\eqref{os3}.
Consider an initial control $v^0$ and assume that, at step $k$, a control
 $v^k$ is known. The computation of $v^{k+1}$ is achieved as follows:
\begin{enumerate}
\item Compute $y^k$, $p^k$ and the associated target trajectory $\chi^k$ according to
Equations \eqref{os1}, \eqref{os2} and \eqref{target_trajectory} respectively.
\item For $n=0,...,N-1$, solve the sub-problems \eqref{min-loc} in
  parallel and denote by $\tilde{v}_n^{k+1}$ the
  corresponding solutions.\label{step2}
\item Define $\tilde v^{k+1}$ as the concatenation of the sequence $(\tilde{v}_n^{k+1})_{n=0,...,N-1}$.
\item\label{step4} Update the control variable by 
\begin{equation}\label{eq:update}
v^{k+1}= v^{k}+\theta_{k}(\tilde
  v^{k+1}-v^{k}),
\end{equation}
where the value $\theta_k$ is chosen to minimizes
  $J(v^{k}+\theta_{k}(\tilde v^{k+1}-v^{k}))$. 
\end{enumerate}
We have not detailed Step \ref{step2} as we rather aim  at
presenting a structure for a general approach. However, because of the strict convexity of the
problems we consider, a small number of conjugate
gradient method steps can be used to achieve the resolution of these
steps.
\section{Numerical experiments}\label{sec:num}
\begin{figure}[!htb]
\psfrag{'./N=1'}[r][c]{\tiny $N=1$ :}
\psfrag{'./N=2'}[r][c]{\tiny$N=2$ :}
\psfrag{'./N=4'}[r][c]{\tiny $N=4$ :}
\psfrag{'./N=8'}[r][c]{\tiny $N=8$ :}
\psfrag{'./N=16'}[r][c]{\tiny $N=16$ :}
\psfrag{'./N=32'}[r][c]{\tiny $N=32$ :}
\psfrag{x}[c]{\footnotesize Number of iterations}
\psfrag{y}[c]{\footnotesize Functional values}
\psfrag{mult}[c]{ \footnotesize Number of multiplications}
\psfrag{T}{\tiny }
\psfrag{Jvalues}[c]{\footnotesize Functional values}
\psfrag{sit}[c]{}
\begin{centering}
\includegraphics[width=.5\textwidth,angle=-90]{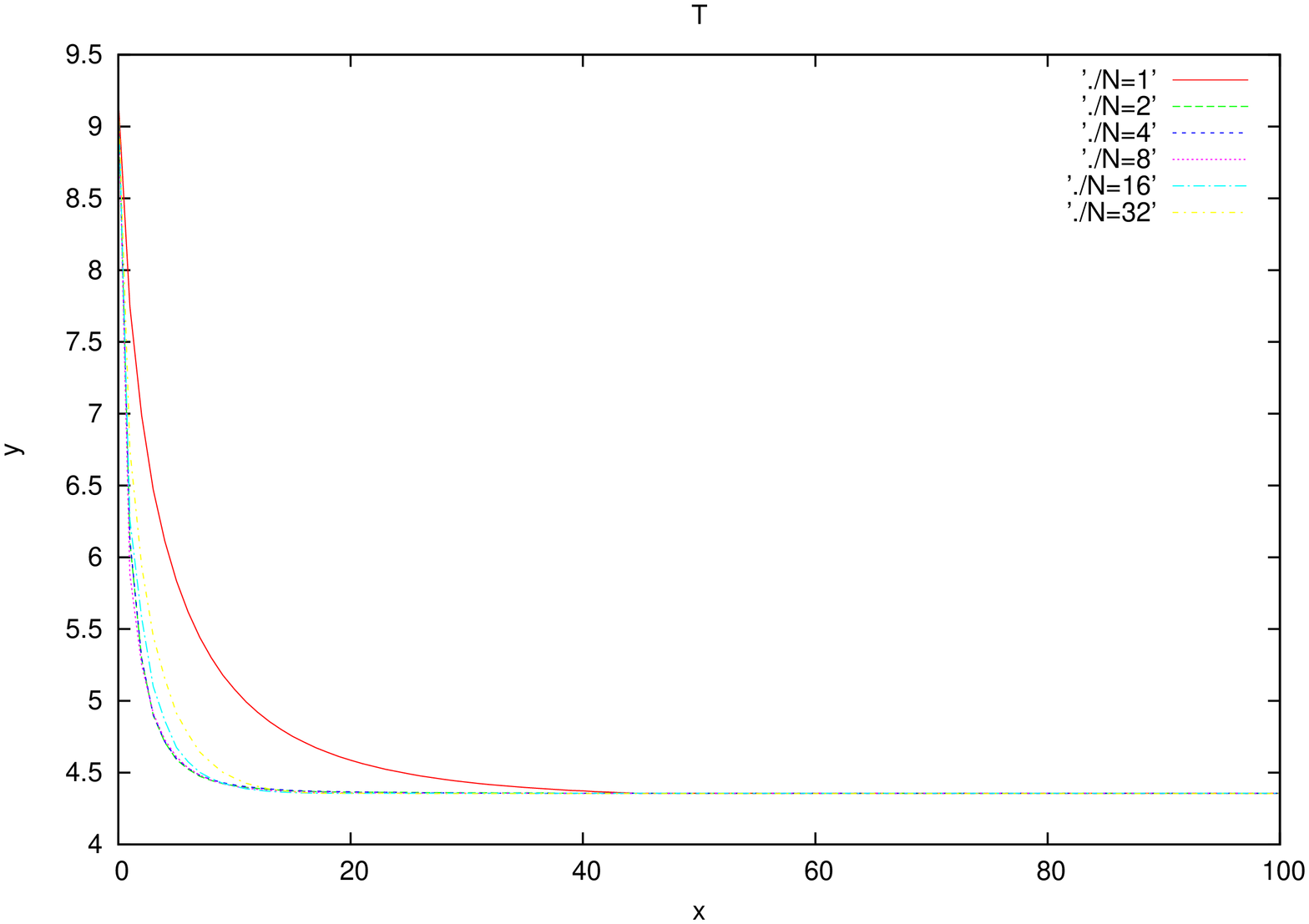}\\
\includegraphics[width=.5\textwidth,angle=-90]{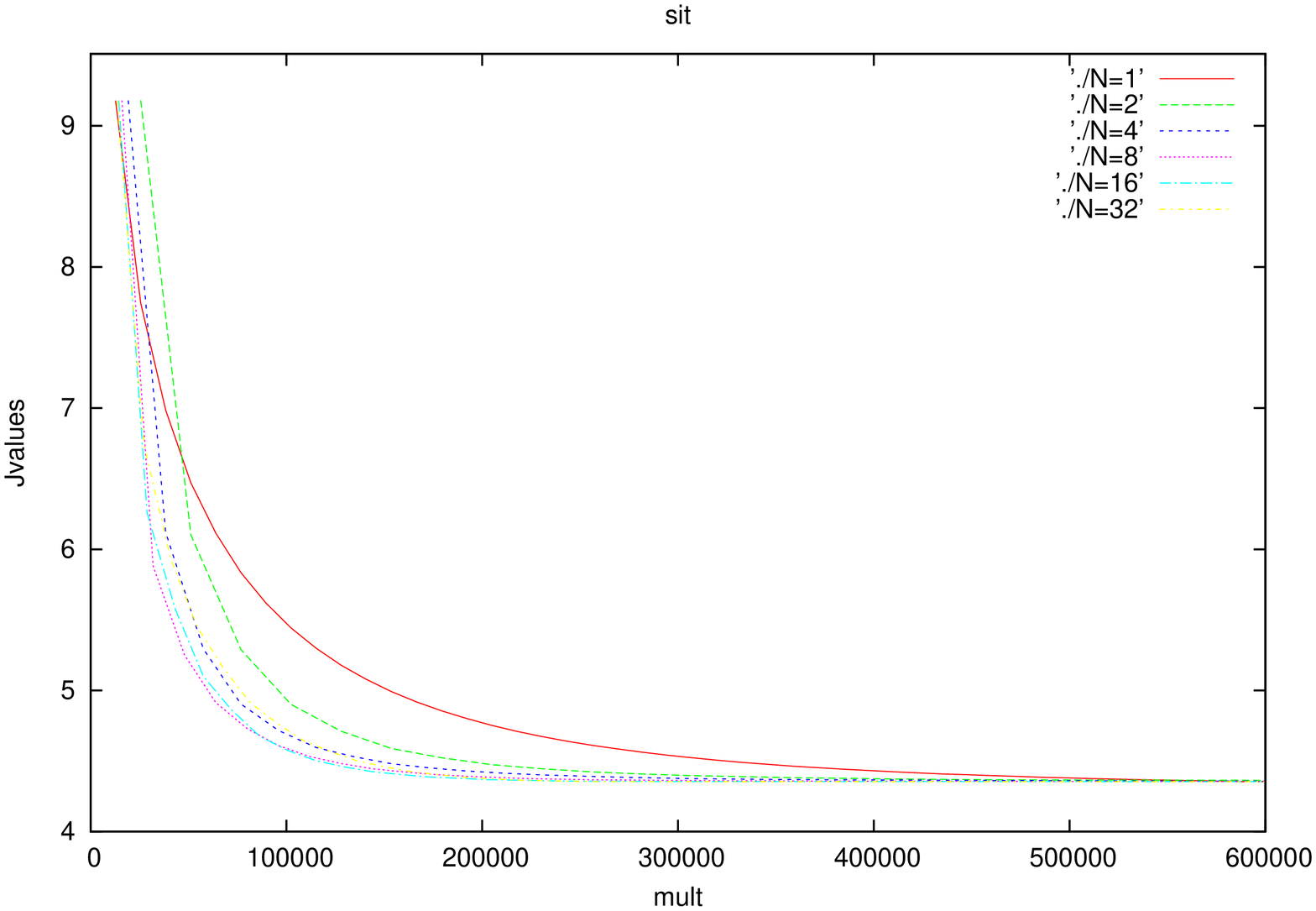}\\
\end{centering}
\caption{Functional values evolution, with respect to the number of
  iterations (top) and multiplications (bottom).}\label{fig}
\end{figure}

In this section, we test the efficiency of our method. We consider a
2D example, where $\Omega=[0,1]\times [0,1]$ and
$\Omega_c=[\frac{1}{3},\frac{2}{3}]\times[\frac{1}{3},\frac{2}{3}]$
. The parameters related to our control problem are $T=6.4$,
$\alpha=10^{-2}$ and $\nu=10^{-2}$. The time interval is discretized using
a uniform step $\delta t=10^{-3}$, and an Implict-Euler solver is used
to approximate the solution of Equations~(\ref{os1}--\ref{os2}). For the space
discretization, we use $\mathbb{P}_{1}$ finite elements. Our
implementation make use of the freeware {\texttt FreeFem}~\cite{bib7}
and the parallelization is achieved thanks to the Message Passing
Interface library. The independent optimization procedures required in
Step~\ref{step2} are simply carried out using one step of an optimal
step gradient method. The results are presented in Figure~\ref{fig}. 

In the first plot, we consider the evolution of the cost functional
values with respect to the iterations and do not take into account the
parallelization. 
The result reveals that our
algorithm significantly accelerates the optimization process. This outcome may
indicates
that the splitting introduced in our approach acts as a preconditionner during the
numerical
optimization. This will be the purpose of some further investigation~\cite{theseK}, in the same spirit as in~\cite{MSS}.

In a second plot, we represent the evolution of the cost functional values
with respect to the number of matrix-vector product. Parallel
computations that are done in Step~~\ref{step2} are only counted
once. When comparing with a standard optimal gradient step method, we
observe speed-up approximatively equal to 3.  


%

\end{document}